%% file: sample.tex
\documentclass[letterpaper, 10 pt, conference]{ieeeconf}  

\IEEEoverridecommandlockouts                              
\title{\LARGE \bf
Smooth and Exact Parameterization of Continuous-time Signal Temporal Logic Specifications for Trajectory Optimization
}

\author{Samet Uzun and Beh\c{c}et A\c{c}\i kme\c{s}e
\thanks{S. Uzun and B. A\c{c}\i kme\c{s}e are with the Department of Aeronautics and Astronautics, University of Washington, Seattle, WA 98105 (email:{samet, behcet}@uw.edu).}%
}

\input{preamble}

\begin{document}

\maketitle
\thispagestyle{empty}
\pagestyle{empty}

\begin{abstract}

This paper presents a smooth parameterization of continuous-time Signal Temporal Logic (CT-STL) specifications for nonconvex trajectory optimization that is sound and complete up to the accuracy of the underlying numerical integration scheme. CT-STL provides a natural framework for encoding rich temporal and logical task requirements, but existing trajectory-optimization formulations typically enforce such specifications only at discrete sampling nodes. In contrast, the proposed method evaluates specifications in dense time, thereby guaranteeing continuous-time satisfaction of \always{} predicates, which is critical for path constraints such as obstacle avoidance, while eliminating the node-induced conservatism of \eventually{} predicates by allowing satisfaction at any time within the prescribed interval. These two dense-time constructions also serve as the main building blocks for handling more general CT-STL formulas, including complex \until{} specifications. Furthermore, the proposed parameterization resolves the locality and gradient-masking issues inherent in standard quantitative semantics, yielding a more favorable landscape for gradient-based solvers. Although dense-time evaluation introduces additional function evaluations during discretization, it also permits substantially coarser temporal grids without sacrificing safety or logical fidelity. This, in turn, reduces the dimension of the resulting nonconvex program, which is often the dominant factor in trajectory-generation cost. The numerical effectiveness and semantic exactness of the proposed framework are demonstrated on an agile quadrotor flight problem subject to a complex continuous-time \until{} specification. The implementation is available at \url{https://github.com/UW-ACL/TrajOpt_CT-STL}

\end{abstract}





\section{INTRODUCTION}

Temporal logic specifications \cite{baier2008principles} provide a systematic framework for encoding high-level tasks and safety requirements in dynamical systems. They have been used in a wide range of autonomy applications, including UAV planning \cite{Pant2018FlyByLogic}, mobile-robot planning \cite{buyukkocak2021planning,cardona2024planning}, and safe learning \cite{aksaray2016q,karagulle2024safe}. 
They have also been incorporated into differentiable and learning-based pipelines for planning and control \cite{leung2023backpropagation}.
While linear temporal logic (LTL) \cite{pnueli1977temporal} and metric interval temporal logic (MITL) \cite{alur1996benefits} provide foundations for specifying temporal requirements over discrete and continuous time, respectively, signal temporal logic (STL) \cite{maler2004monitoring} is particularly well suited to physical systems because its predicates are defined directly over real-valued continuous signals. A key advance in STL-based control synthesis was the introduction of quantitative semantics, or robustness \cite{donze2010robust}. Using $\min$ and $\max$ operations, robustness assigns a continuous measure of satisfaction or violation, which enables temporal specifications to be incorporated into numerical optimization as constraints or objective terms \cite{raman2014model,uzun2024optimization}.

To solve these optimization problems, the vast majority of existing literature relies on discrete-time evaluations of STL (DT-STL). The exact $\min$ and $\max$ functions of discrete-time space robustness (D-SR) are traditionally encoded using mixed-integer programming (MIP) \cite{raman2014model,sadraddini2015robust,rodionova2021time}. However, the worst-case computational complexity of MIP formulations grows exponentially with the number of binary variables. Alternatively, D-SR can be used directly in continuous optimization, but the inherent nonsmoothness of the operators prevents the use of gradient-based optimization algorithms. Various smooth approximations, such as log-sum-exponential (LSE) \cite{pant2017smooth}, and softmax \cite{gilpin2020smooth}, restore gradient information but introduce approximation errors and frequently suffer from \textit{locality} and \textit{masking} issues \cite{mehdipour2019average}, where gradients vanish and stall solver convergence. 

Recently, the discrete-time generalized mean-based smooth robustness (D-GMSR) \cite{uzun2024optimization} was introduced to definitively resolve these discrete-time trade-offs. By composing specialized generalized mean functions, D-GMSR provides an exact, $\mathcal{C}^1$-smooth parameterization of STL robustness that eliminates approximation errors while naturally distributing gradient information to overcome locality and masking. This formulation enables the use of high-performance sequential convex programming (SCP) solvers \cite{proxconvex} for highly complex DT-STL trajectory optimization problems.

Despite advancements in discrete-time smooth formulations, evaluating STL specifications solely at sparse discrete nodes leaves the continuous-time inter-sample behavior of the system unchecked. Beyond discrete-time encodings, a growing body of literature addresses dense-time satisfaction for continuous-time systems. In control synthesis, continuous-time STL satisfaction is frequently converted into time-varying set-invariance or reachability conditions enforced via control barrier functions (CBFs) \cite{Lindemann2019CBFSTL,Charitidou2021BFMPCSTL,Buyukkocak2022ActuationCBF,Yu2024NestedSTL}. However, these approaches focus on online controller synthesis rather than providing a direct, smooth parameterization of STL constraints for trajectory optimization.

For trajectory optimization, continuous-time STL guarantees have historically relied on specialized, often restrictive mechanisms. Some approaches augment mixed-integer planners with CBF-derived intersample constraints under a zero-order hold \cite{Yang2020CTSTLCBF}, limiting scalability to linear systems and predicates. Others plan over sparse waypoints and impose strictified sampled-time specifications alongside hierarchical tracking controllers \cite{Pant2018FlyByLogic}, introducing conservatism and platform-specific assumptions. Recent works optimizing directly over continuous parameterizations employ B\'ezier curves with mixed-integer linear programs (MILPs) \cite{Verhagen2024TemporallyRobust} or time-varying robustness constraints \cite{Yuan2024TimeVaryingRobustness,Yuan2025Quadrotors}, or utilize sampling-based planners searching within forward-invariant sets \cite{Marchesini2025SamplingBased}. While advancing continuous-time planning, these guarantees remain tightly coupled to specific machinery, such as barrier certificates, tracking-error bounds, piecewise polynomials, or mixed-integer encodings, falling short of a general, smooth, and semantically exact parameterization. 

Building on this line of work, a continuous-time successive convexification framework for enforcing the continuous-time satisfaction of path constraints, i.e., \always{}-type specifications, was developed in~\cite{ctcs2024}. Their numerical effectiveness has been demonstrated across a range of applications, including GPU-accelerated trajectory optimization for 6-DoF rocket landing \cite{chari2024fast}, nonlinear model predictive control (NMPC) for obstacle avoidance \cite{nmpc2024}, and trajectory optimization for 6-DoF aircraft approach and landing \cite{aircraft2024}. This framework was subsequently combined with GMSR in~\cite{uzun2025sequential} to handle the continuous-time \always{} satisfaction of logical specifications in a 6-DoF rocket landing problem. Its practical utility was further illustrated in perception-constrained drone motion planning \cite{uzun2025motion}, 
a rocket-landing application solved via an auto-tuned primal-dual algorithm~\cite{mceowen2026auto}, and SCvxGEN, an automatic code-generation framework for fast embedded trajectory optimization \cite{scvxgen}. Despite these advances, existing GMSR-based continuous-time formulations remain limited primarily to \always{}-type specifications, and a general treatment of richer temporal operators such as \eventually{} and \until{} is still lacking.

Motivated by this gap, this paper develops a smooth parameterization of continuous-time Signal Temporal Logic (CT-STL) specifications for nonconvex trajectory optimization. The proposed framework extends the exact and smooth D-GMSR formulation \cite{uzun2024optimization} from discrete-time to dense-time STL semantics, thereby enabling continuous-time satisfaction of rich temporal and logical task requirements within a standard nonconvex optimization framework. In contrast to existing CT-STL approaches that rely on conservative inter-sample bounds, specialized trajectory parameterizations, or mixed-integer machinery, the proposed method evaluates STL specifications directly on a dense-time representation of the trajectory. As a result, it guarantees continuous-time satisfaction of \always{}-type specifications, which is critical for safety constraints such as obstacle avoidance, while removing the node-induced conservatism of \eventually{}-type specifications by allowing satisfaction at any time within the prescribed interval. Moreover, by inheriting the key properties of GMSR, the resulting formulation remains free from the locality and masking issues that commonly degrade standard quantitative semantics, yielding a more favorable landscape for gradient-based solvers. Owing to the exactness properties of GMSR, the proposed CT-STL parameterization is sound and complete with respect to the underlying dense-time semantics up to the numerical accuracy of the integration scheme used to represent the continuous-time dynamics.

\section{Preliminaries} \label{sec:stl}

This section briefly reviews the syntax and semantics of Signal Temporal Logic (STL). Let $x : \mathbb{R}_{+} \to \mathbb{R}^n$ be a continuous-time signal. A Boolean-valued predicate $\mu$ is induced by a real-valued $\mathcal{C}^1$ predicate function $g:\mathbb{R}^n \to \mathbb{R}$ according to
\begin{align*}
    \mu := (g(x(t)) \ge 0).
\end{align*}
That is, $\mu$ evaluates to \stltrue{} if $g(x(t))\ge 0$ and \stlfalse{} otherwise.

The syntax of an STL formula $\varphi$ is defined recursively as
\begin{equation*} \label{eq:stl_syntax}
    \varphi ::= \top \mid \mu \mid \neg \varphi \mid \varphi_1 \wedge \varphi_2 \mid \varphi_1 \bm{U}_{[a,b]} \varphi_2,
\end{equation*}
where $\varphi_1$ and $\varphi_2$ are STL subformulas, $\neg$ and $\wedge$ are the fundamental Boolean operators, and $\bm{U}_{[a,b]}$ is the bounded \until{} operator over $[a,b]\subset\mathbb{R}_+$.

The standard derived operators \stlor{}, \implication{}, \eventually{}, and \always{} are obtained from the fundamental operators through
\[
\varphi_1 \vee \varphi_2 \equiv \neg(\neg \varphi_1 \wedge \neg \varphi_2),
\qquad
\varphi_1 \implies \varphi_2 \equiv \neg \varphi_1 \vee \varphi_2,
\]
and
\[
\bm{F}_{[a,b]}\varphi \equiv \top \bm{U}_{[a,b]} \varphi,
\qquad
\bm{G}_{[a,b]}\varphi \equiv \neg \bm{F}_{[a,b]} \neg \varphi.
\]

The Boolean semantics specify when a signal $x$ satisfies a formula $\varphi$ at time $t$, denoted by $(x,t)\models \varphi$.

\begin{definition}[Continuous-time STL semantics \cite{donze2010robust}]
For a continuous-time signal $x:\mathbb{R}_{+}\to\mathbb{R}^n$ and time $t\in\mathbb{R}_+$, the Boolean satisfaction relation $(x,t)\models\varphi$ is defined recursively as in Table~\ref{tab:stl_semantics}.
\end{definition}

\begin{table}[H]
\centering
\renewcommand{\arraystretch}{1.25}
\begin{tabular}{ll}
\toprule
\textbf{Formula} $\varphi$ & \textbf{Boolean semantics} $(x,t)\models \varphi$ \\
\midrule
$\mu$ & $g(x(t)) \ge 0$ \\
$\neg \varphi$ & $\neg\big((x,t)\models \varphi\big)$ \\
$\varphi_1 \wedge \varphi_2$ & $(x,t)\models \varphi_1 \;\wedge\; (x,t)\models \varphi_2$ \\
$\varphi_1 \vee \varphi_2$ & $(x,t)\models \varphi_1 \;\vee\; (x,t)\models \varphi_2$ \\
$\varphi_1 \implies \varphi_2$ & $\neg\big((x,t)\models \varphi_1\big)\;\vee\; (x,t)\models \varphi_2$ \\
\midrule
$\bm{G}_{[a,b]}\varphi$ & $\forall\, t' \in [t+a,t+b],\; (x,t')\models \varphi$ \\
$\bm{F}_{[a,b]}\varphi$ & $\exists\, t' \in [t+a,t+b]\ \text{s.t.}\ (x,t')\models \varphi$ \\
$\varphi_1 \bm{U}_{[a,b]} \varphi_2$ &
$\exists\, t_1 \in [t+a,t+b]\ \text{s.t.}\ (x,t_1)\models \varphi_2$ \\
& \quad and $\forall\, t_2 \in [t,t_1],\; (x,t_2)\models \varphi_1$ \\
\bottomrule
\end{tabular}
\caption{Continuous-time Boolean semantics for STL.}
\label{tab:stl_semantics}
\end{table}

While Boolean semantics provide a \stltrue{}/\stlfalse{} evaluation, optimal control requires a continuous-valued measure that can guide gradient-based optimization. A standard choice is the quantitative robustness semantics introduced in \cite{donze2010robust}. Let $\rho^\varphi(x,t)$ denote the robustness of formula $\varphi$ at time $t$. The sign of $\rho^\varphi(x,t)$ is consistent with Boolean satisfaction, and in particular $\rho^\varphi(x,t) > 0$ implies strict satisfaction. The standard robustness semantics are summarized in Table~\ref{tab:stl_robustness}.

\begin{table}[H]
\centering
\renewcommand{\arraystretch}{1.25}
\begin{tabular}{ll}
\toprule
\textbf{Formula} $\varphi$ & \textbf{Quantitative robustness} $\rho^\varphi(x,t)$ \\
\midrule
$\mu$ & $g(x(t))$ \\
$\neg \varphi$ & $-\rho^\varphi(x,t)$ \\
$\varphi_1 \wedge \varphi_2$ & $\min\!\big(\rho^{\varphi_1}(x,t),\,\rho^{\varphi_2}(x,t)\big)$ \\
$\varphi_1 \vee \varphi_2$ & $\max\!\big(\rho^{\varphi_1}(x,t),\,\rho^{\varphi_2}(x,t)\big)$ \\
$\varphi_1 \implies \varphi_2$ & $\max\!\big(-\rho^{\varphi_1}(x,t),\,\rho^{\varphi_2}(x,t)\big)$ \\
\midrule
$\bm{G}_{[a,b]}\varphi$ & $\min\limits_{t' \in [t+a,t+b]} \rho^\varphi(x,t')$ \\
$\bm{F}_{[a,b]}\varphi$ & $\max\limits_{t' \in [t+a,t+b]} \rho^\varphi(x,t')$ \\
$\varphi_1 \bm{U}_{[a,b]} \varphi_2$
& $\max\limits_{t_1 \in [t+a,t+b]}
\min\!\Big(\rho^{\varphi_2}(x,t_1),\,
\min\limits_{t_2 \in [t,t_1]} \rho^{\varphi_1}(x,t_2)\Big)$ \\
\bottomrule
\end{tabular}
\caption{Standard quantitative robustness semantics for continuous-time STL.}
\label{tab:stl_robustness}
\end{table}

\section{Trajectory Optimization with CT-STL Specifications}

We consider a fixed-final-time optimal control problem in Mayer form, subject to continuous-time Signal Temporal Logic (CT-STL) specifications:
\begin{subequations} \label{ct-ocp}
    \begin{align}
        \underset{x,\,u}{\mathrm{minimize}} \quad
        & L\big(x(t_{\mathrm f})\big) \label{ct-cost} \\
        \mathrm{subject~to} \quad
        & \dot{x}(t)=F\big(t,x(t),u(t)\big),
        \; \text{a.e. } t\in[0,t_{\mathrm f}], \\
        & P\big(x(0),u(0),x(t_{\mathrm f}),u(t_{\mathrm f})\big)\le 0, \\
        & Q\big(x(0),u(0),x(t_{\mathrm f}),u(t_{\mathrm f})\big)=0, \\
        & (x,0)\models \varphi. \label{eq:ct-stl-constraint}
    \end{align}
\end{subequations}
Here, 
$L$ denotes the terminal cost, $F$ defines the nonlinear system dynamics, $P$ and $Q$ represent the boundary inequality and equality constraints, and $\varphi$ is the STL specification evaluated at the initial time.

\begin{remark} \label{rem:ctcs-always}
Continuous-time path constraints can be encoded directly as STL formulas. For example, an inequality constraint $g(x(t))\le 0$ and an equality constraint $h(x(t))=0$ enforced over the entire trajectory can be written as
\begin{align*}
    \bm G_{[0,t_{\mathrm f}]}
    \Big(
        (-g(x)\ge 0)\wedge (h(x)\ge 0)\wedge (-h(x)\ge 0)
    \Big).
\end{align*}
Thus, path constraints can be handled within the same smooth robustness framework as the STL task specifications.
\end{remark}

\subsection{Control Input Parameterization}\label{ssec:control_param}

To obtain a finite-dimensional optimization problem, we parameterize the continuous-time control input using a finite set of nodal values. Specifically, we discretize the fixed time horizon $[0,t_{\mathrm f}]$ using $K$ nodes
\begin{align*}
    0=t_1<t_2<\cdots<t_K=t_{\mathrm f},
\end{align*}
with uniform spacing
\begin{align*}
    \Delta t := \frac{t_{\mathrm f}}{K-1},
    \qquad
    t_k=(k-1)\Delta t,
    \quad k=1,\dots,K.
\end{align*}
The control input is then parameterized by first-order hold (FOH), so that for each interval $[t_k,t_{k+1}]$,
\begin{align*}
    u(t)
    =
    u_k \frac{t_{k+1}-t}{\Delta t}
    +
    u_{k+1}\frac{t-t_k}{\Delta t},
    \qquad t\in[t_k,t_{k+1}],
\end{align*}
where $u_k\in\mathbb{R}^m$ denotes the control value at node $t_k$, and $k\in\mathcal K:=\{1,\dots,K-1\}$. This parameterization replaces the infinite-dimensional control trajectory $u(\cdot)$ with the finite set of decision variables $\{u_k\}_{k=1}^K$, while retaining a continuous-time representation of the control signal between nodes.

\subsection{Discretization and Multiple Shooting}\label{ssec:discretization}

With the state and control values $x_k$ and $u_k$ at the discrete nodes $t_k$ acting as decision variables, we discretize the system dynamics using a multiple-shooting method~\cite{bock1984multiple}. The state transition from node $k$ to node $k+1$ is enforced through
\begin{align*}
    x_{k+1} = \Phi_{t_k}^{t_{k+1}}(x_k,u_k,u_{k+1}),
\end{align*}
where the flow map is defined by
\begin{align*}
    \Phi_{t_k}^{t_{k+1}}(x_k,u_k,u_{k+1})
    :=
    x_k + \int_{t_k}^{t_{k+1}} F\big(t,x(t),u(t)\big)\,\mathrm dt.
\end{align*}
Here, $x(t)$ is the continuous-time state trajectory satisfying $x(t_k)=x_k$ under the FOH control parameterization over $[t_k,t_{k+1}]$.

To evaluate this transition numerically, we partition each shooting interval $[t_k,t_{k+1}]$ into $N$ integration substeps. Let the corresponding integration subnodes be
\begin{align}
    t_k=t_k^0<t_k^1<\cdots<t_k^N=t_{k+1}.
    \label{eq:disc_t}
\end{align}
The state at each intermediate subnode $t_k^i$ is obtained by applying a chosen numerical integration scheme over these substeps, for example a Runge--Kutta method. Thus, each subnode state $x(t_k^i)$ is generated explicitly from $x_k$ under the FOH control parameterization and is therefore a well-defined, continuously differentiable function of the nodal variables, provided the underlying integration scheme is itself differentiable. Consequently, its Jacobians with respect to $x_k$, $u_k$, and $u_{k+1}$ are well defined and can be evaluated efficiently, which is important for gradient-based optimization. The terminal subnode value $x(t_k^N)$ provides the numerical evaluation of the flow map and is used to enforce the multiple-shooting defect constraint against the subsequent decision variable $x_{k+1}$.

\subsection{Parameterization of the CT-STL specifications via GMSR}

To incorporate CT-STL specifications into trajectory optimization, we use the generalized mean-based smooth robustness (GMSR) construction introduced in \cite{uzun2024optimization}. Standard quantitative semantics typically face a fundamental trade-off. They are either exact but nonsmooth, due to repeated use of $\min$ and $\max$, or smooth but approximate. Moreover, common smooth approximations often suffer from \emph{locality}, where sensitivity is concentrated at a single critical time, and \emph{masking}, where satisfied subformulas suppress gradient information from violated ones. GMSR avoids these issues by providing a smooth and exact parameterization of the logical operators, with gradient information distributed across multiple relevant subformulas and time instants rather than collapsing onto a single critical sample. When evaluated on the dense numerical trajectory induced by the dynamics discretization, it yields an exact CT-STL realization on that trajectory representation.

\subsubsection{Parameterization of the logical operators}

For the logical operators $\wedge$ and $\vee$, we use the following simplified GMSR parameterization:
\begin{align*}
    {}^{\wedge} h^{c}(y)
    &\defeq
    \left( M_{0}^{c}\!\left(\relu{y}^{2}\right) \right)^{\frac12}
    -
    \left( M_{1}^{c}\!\left(\negrelu{y}^{2}\right) \right)^{\frac12}, \\
    {}^{\vee} h^{c}(y)
    &\defeq
    -{}^{\wedge} h^{c}(-y),
\end{align*}
where $y\in\mathbb{R}^n$, $\relu{y}:=\max(0,y)$, and $\negrelu{y}:=\min(0,y)$ are applied elementwise, and
\begin{align*}
    M_{0}^{c}(z)
    &\defeq
    \left(c^n+\prod_{i=1}^n z_i\right)^{1/n}, \\
    M_{1}^{c}(z)
    &\defeq
    c+\frac{1}{n}\sum_{i=1}^n z_i,
\end{align*}
with $c\in\mathbb{R}_{++}$.

\begin{remark}[Smoothing parameter and GMSR extensions]
The parameter $c$ acts as a positive regularization shift near the switching region. For any $c>0$, the resulting functions ${}^{\wedge} h^{c}$ and ${}^{\vee} h^{c}$ are $\mathcal{C}^1$-smooth with bounded gradients, while preserving the exact sign-based parameterization of \conjunction{} and \disjunction{}. Smaller values of $c$ produce a less flattened transition near zero, whereas larger values yield a smoother but flatter local landscape. 

The simplified form above corresponds to a uniform-weight, fixed-curvature choice within the more general GMSR family introduced in \cite{uzun2024optimization}. In that formulation, weight parameters can be used to emphasize selected inputs or time samples, while curvature parameters can make the aggregation more or less extremum-like. 
For clarity of presentation, we use the simplified form throughout this paper. Further intuition, derivative properties, and plots for different choices of these parameters are given in \cite{uzun2024optimization}.
\end{remark}

We now introduce the robustness measure recursively. For an atomic predicate $\mu := (g(x)\ge 0)$, define
\begin{align*}
    \Gamma^{\mu}(x,t) &:= g(x(t)).
\end{align*}
\Negation{} ($\neg$), \conjunction{} ($\wedge$), \disjunction{} ($\vee$), and \implication{} ($\implies$) are then parameterized as
\begin{align*}
    \Gamma_{\bm c}^{\neg\varphi}(x,t)
    &:= -\Gamma_{\bm c}^{\varphi}(x,t), \\
    \Gamma_{\bm c}^{\varphi_1 \wedge \varphi_2}(x,t)
    &:= {}^{\wedge} h^{c_1}
    \Big(
        \Gamma_{\bm c}^{\varphi_1}(x,t),\,
        \Gamma_{\bm c}^{\varphi_2}(x,t)
    \Big), \\
    \Gamma_{\bm c}^{\varphi_1 \vee \varphi_2}(x,t)
    &:= {}^{\vee} h^{c_1}
    \Big(
        \Gamma_{\bm c}^{\varphi_1}(x,t),\,
        \Gamma_{\bm c}^{\varphi_2}(x,t)
    \Big), \\
    \Gamma_{\bm c}^{\varphi_1 \Rightarrow \varphi_2}(x,t)
    &:= {}^{\vee} h^{c_1}
    \Big(
        -\Gamma_{\bm c}^{\varphi_1}(x,t),\,
        \Gamma_{\bm c}^{\varphi_2}(x,t)
    \Big),
\end{align*}
where $\bm c$ collects the positive shift parameters used in the recursive GMSR compositions.

Accordingly, for predicates of the form $\varphi_i := (g_i(x)\ge 0)$, satisfaction at time $t$ is characterized exactly by
\begin{align*}
    (x,t)\models \bigwedge_{i=1}^{n}\varphi_i
    &\iff
    {}^{\wedge} h^{c}\!\left(
        g_1(x(t)),\dots,g_n(x(t))
    \right)\ge 0, \\
    (x,t)\models \bigvee_{i=1}^{n}\varphi_i
    &\iff
    {}^{\vee} h^{c}\!\left(
        g_1(x(t)),\dots,g_n(x(t))
    \right)\ge 0.
\end{align*}
More generally,
\begin{align*}
    (x,t)\models \varphi
    \iff
    \Gamma_{\bm c}^{\varphi}(x,t)\ge 0
\end{align*}
for any formula $\varphi$ built recursively from the above logical operators.

\begin{remark}
Although the conditions
\begin{align*}
    {}^{\wedge} h^{c}(y)\ge 0
    &\iff
    \sum_{i=1}^{n}\negrelu{y_i}^{2}=0, \\
    {}^{\vee} h^{c}(y)\ge 0
    &\iff
    \prod_{i=1}^{n}\negrelu{y_i}^{2}=0
\end{align*}
are logically equivalent, the full GMSR expressions are numerically preferable. The geometric-mean normalization, together with the outer square-root, keeps the robustness values and gradients better scaled, while the equivalent sum and product tests can be poorly conditioned. Moreover, the shift parameter $c>0$ regularizes the switching region, ensures $\mathcal{C}^1$-smoothness, and avoids singular behavior near zero.
\end{remark}

\subsubsection{Parameterization of the temporal operators}

We parameterize the temporal operators directly on the dense-time grid induced by the numerical integration of the dynamics. In direct optimal control, continuous-time satisfaction between mesh points is often enforced through conservative analytical bounds, which shrink the feasible set and may degrade optimality. Here, instead, we evaluate the STL semantics directly on the dense numerical mesh already used to integrate the dynamics. This yields an exact CT-STL realization on the chosen dense-time trajectory representation, without introducing any additional approximation parameter beyond the discretization of the physical dynamics themselves.

To avoid introducing a separate discretization solely for STL evaluation, we reuse the intermediate integration subnodes from \eqref{eq:disc_t}. Let these subnodes be flattened into a single ordered sequence
\[
\{\bar t_m\}_{m=1}^{M},
\qquad
M=(K-1)N+1,
\]
and let
$
\bar x_m := x(\bar t_m)
$
denote the corresponding states. 
The following schematic illustrates how the nodal states \(x_k\), the intermediate subnodes \(t_k^i\), and the flattened dense-time sequence \(\bar{x}_m=x(\bar t_m)\) are related.
\vspace{-0.4cm}
\begin{center}
\begin{tikzpicture}[x=1.0cm,y=0.7cm,>=latex]

\definecolor{mainblue}{RGB}{40,90,160}
\definecolor{suborange}{RGB}{220,130,40}
\definecolor{densegreen}{RGB}{20,120,80}

\draw[thick] (0,0) -- (8,0);

\foreach \x in {0,1.0,2.0,2.7,3.8,4.8,5.8,6.5,8.0}
    \fill (\x,0) circle (1.2pt);

\node[above=3pt, text=mainblue] at (0,0) {$x_1$};
\node[above=3pt, text=mainblue] at (3.8,0) {$x_2$};
\node[above=3pt, text=mainblue] at (8.0,0) {$x_K$};

\node[above=3pt, text=suborange] at (1.0,0) {$t_1^1$};
\node[above=3pt, text=suborange] at (2.0,0) {$t_1^2$};
\node[above=3pt] at (2.7,0) {$\cdots$};
\node[above=3pt, text=suborange] at (4.8,0) {$t_2^1$};
\node[above=3pt, text=suborange] at (5.8,0) {$t_2^2$};
\node[above=3pt] at (6.5,0) {$\cdots$};

\node[below=3pt, text=densegreen] at (0,0) {$\bar{x}_1$};
\node[below=3pt, text=densegreen] at (1.0,0) {$\bar{x}_2$};
\node[below=3pt, text=densegreen] at (2.0,0) {$\bar{x}_3$};
\node[below=3pt] at (2.7,0) {$\cdots$};
\node[below=3pt, text=densegreen] at (3.8,0) {$\bar{x}_{N+1}$};
\node[below=3pt, text=densegreen] at (4.8,0) {$\bar{x}_{N+2}$};
\node[below=3pt, text=densegreen] at (5.8,0) {$\bar{x}_{N+3}$};
\node[below=3pt] at (6.5,0) {$\cdots$};
\node[below=3pt, text=densegreen] at (8.0,0) {$\bar{x}_M$};

\end{tikzpicture}
\end{center}
\vspace{-0.2cm}
Because each dense-time state $\bar x_m$ is generated by numerical integration from the nodal decision variables, the resulting robustness values are smooth functions of the discretized trajectory variables.

Let
\[
X := \{x_k\}_{k=1}^{K},
\qquad
U := \{u_k\}_{k=1}^{K},
\]
and write
$
\Gamma_{\bm c}^{\varphi}(X,U;\bar t_m)
$
for the robustness value of subformula $\varphi$ at dense-time node $\bar t_m$, where the dependence on $(X,U)$ enters through the dense numerical trajectory. For an atomic predicate $\mu := (g(x)\ge 0)$,
\[
\Gamma^{\mu}(X,U;\bar t_m) := g(\bar x_m).
\]
For composite formulas, the values $\Gamma_{\bm c}^{\varphi}(X,U;\bar t_m)$ are obtained recursively from the logical constructions above.

Now consider a formula evaluated at dense-time node $\bar t_\ell$. For a time interval $[a,b]$, define
\[
\mathcal I_\ell[a,b]
:=
\left\{
m\in\{1,\dots,M\}
:\;
\bar t_m \in [\bar t_\ell+a,\ \bar t_\ell+b]
\right\},
\]
assuming the interval lies within the time horizon. If the boundary points are not already included in the dense-time grid, they can be inserted by refining the corresponding integration subintervals.

Using this notation, the temporal robustness measures are defined by
\begin{align*}
    \Gamma_{\bm c}^{\bm F_{[a,b]}\varphi}(X,U;\bar t_\ell)
    &:=
    {}^{\vee} h^{c_1}
    \Big(
        [\,\Gamma_{\bm c}^{\varphi}(X,U;\bar t_m)\,]_{m\in\mathcal I_\ell[a,b]}
    \Big), \\
    \Gamma_{\bm c}^{\bm G_{[a,b]}\varphi}(X,U;\bar t_\ell)
    &:=
    {}^{\wedge} h^{c_1}
    \Big(
        [\,\Gamma_{\bm c}^{\varphi}(X,U;\bar t_m)\,]_{m\in\mathcal I_\ell[a,b]}
    \Big).
\end{align*}
Accordingly,
\begin{align*}
    (x,\bar t_\ell)\models \bm F_{[a,b]}\varphi
    &\iff
    \Gamma_{\bm c}^{\bm F_{[a,b]}\varphi}(X,U;\bar t_\ell)\ge 0, \\
    (x,\bar t_\ell)\models \bm G_{[a,b]}\varphi
    &\iff
    \Gamma_{\bm c}^{\bm G_{[a,b]}\varphi}(X,U;\bar t_\ell)\ge 0.
\end{align*}

Similarly, the \until{} operator is parameterized as
\begin{align*}
    \Gamma_{\bm c}^{\varphi_1 \bm U_{[a,b]} \varphi_2}(X,U;\bar t_\ell)
    &:=
    {}^{\vee} h^{c_1}
    \Big(
        [\,z_m\,]_{m\in\mathcal I_\ell[a,b]}
    \Big),
\end{align*}
where
\begin{align*}
    z_m
    &:=
    {}^{\wedge} h^{c_2}
    \Big(
        \Gamma_{\bm c}^{\varphi_2}(X,U;\bar t_m),\;
        {}^{\wedge} h^{c_3}
        \big(
            [\,\Gamma_{\bm c}^{\varphi_1}(X,U;\bar t_q)\,]_{q=\ell}^{m}
        \big)
    \Big).
\end{align*}
Hence,
\begin{align*}
    (x,\bar t_\ell)\models \varphi_1 \bm U_{[a,b]} \varphi_2
    &\iff
    \Gamma_{\bm c}^{\varphi_1 \bm U_{[a,b]} \varphi_2}(X,U;\bar t_\ell)\ge 0.
\end{align*}

Thus, $\bm G$ requires satisfaction at all dense-time samples in the interval, $\bm F$ requires satisfaction at at least one such sample, and $\bm U$ requires the existence of a sample at which $\varphi_2$ holds while $\varphi_1$ holds at all preceding samples.

The resulting parameterization is smooth with respect to the optimization variables, since the dependence on the decision variables enters only through the integrated dense-time states and the smooth GMSR operators. Moreover, because the temporal operators are evaluated on the same dense-time grid used to discretize the dynamics, the CT-STL semantics are enforced exactly on the dense numerical trajectory representation, up to the numerical accuracy of the underlying integration scheme.

\subsubsection{Gradient behavior and optimization relevance}

A key optimization advantage of GMSR is that its gradient is distributed across multiple relevant arguments, rather than collapsing to a single min/max witness. This property holds both for logical operators, where the inputs are subformula robustness values, and for temporal operators, where the inputs are robustness values sampled along the dense-time trajectory. As a result, the same gradient-distribution mechanism appears across subformulas in the logical case and across time samples in the temporal case, which mitigates locality and masking.

For \conjunction{}, the positive branch is governed by the geometric term and the negative branch by the arithmetic-mean term. Up to a common positive scaling factor, the gradient magnitudes with respect to the inputs satisfy
\begin{align*}
    \left|\frac{\partial\, {}^{\wedge} h^{c}(y)}{\partial y_i}\right|
    &\propto \frac{1}{y_i},
    \qquad y_i > 0, \\
    \left|\frac{\partial\, {}^{\wedge} h^{c}(y)}{\partial y_i}\right|
    &\propto |y_i|,
    \qquad y_i < 0.
\end{align*}
Hence, when a \conjunction{} is satisfied, smaller positive margins receive larger emphasis; when it is violated, larger violations receive larger corrective push. In other words, $\wedge$ emphasizes the weakest satisfied input and, if infeasible, pushes the worst violating input more strongly. Importantly, however, all relevant inputs still receive gradient information.

For \disjunction{}, the behavior is dual:
\begin{align*}
    \left|\frac{\partial\, {}^{\vee} h^{c}(y)}{\partial y_i}\right|
    &\propto |y_i|,
    \qquad y_i > 0, \\
    \left|\frac{\partial\, {}^{\vee} h^{c}(y)}{\partial y_i}\right|
    &\propto \frac{1}{|y_i|},
    \qquad y_i < 0.
\end{align*}
Thus, when a \disjunction{} is violated, the inputs closest to satisfaction receive the largest emphasis, while the remaining ones still contribute. When it is satisfied, stronger satisfying witnesses receive larger emphasis. Accordingly, $\vee$ tends to improve the best candidate first, but without restricting the gradient to a single input.

The temporal operators inherit exactly the same structure, since they are obtained by applying these logical aggregations to the vector of robustness values over the dense-time samples in the relevant interval. Therefore, for $\bm G$, gradient information is distributed over all sampled times in the interval, with stronger emphasis on the smallest positive margins when the specification is satisfied and on the worst violations when it is not. For $\bm F$, the gradient is again distributed over all sampled times, but now the most promising sample receives the largest emphasis. If the specification is violated, this is the sample closest to satisfaction; if it is satisfied, this is the strongest witness sample.

In summary, the same mechanism applies in both dimensions of the construction. Across subformulas for logical operators and across time for temporal operators. An \always{}-type condition pushes the worst parts more strongly, whereas an \eventually{}-type condition pushes the best candidates more strongly, while still propagating gradient information to all relevant components. This is precisely why GMSR alleviates locality and masking and is well-suited for gradient-based trajectory optimization.

\section{Numerical Example}\label{sec:numerical}

We demonstrate the proposed framework on a quadrotor guidance problem with a charging-station waypoint. To keep the example simple, we use a point-mass model of the quadrotor and consider a feasibility problem with fixed initial and terminal states. The task is to find a dynamically feasible trajectory that reaches the terminal state while satisfying a continuous-time \until{} specification, namely that the vehicle speed must remain below a prescribed threshold until the quadrotor reaches the charging station. In addition, continuous-time bounds on tilt angle, thrust magnitude, and speed are enforced through \always{} operators. The implementation is available at \url{https://github.com/UW-ACL/TrajOpt_CT-STL}

\subsection{Problem setup}\label{sec:num_setup}

We use a point-mass quadrotor model with state and control
\begin{align*}
    x := (r,v) \in \mathbb{R}^6,
    \qquad
    u := (u_x,u_y,u_z) \in \mathbb{R}^3,
\end{align*}
where $r=(r_x,r_y,r_z)$ is the position, $v=(v_x,v_y,v_z)$ is the velocity, and $u$ is the thrust vector. The dynamics are
\begin{align}
    \dot{x}(t)
    =
    \begin{bmatrix}
        \dot r(t) \\
        \dot v(t)
    \end{bmatrix}
    =
    \begin{bmatrix}
        v(t) \\
        \frac{1}{m}u(t) + g
    \end{bmatrix},
    \label{eq:num_double_integrator}
\end{align}
where $m>0$ is the vehicle mass and $g=(0,0,-g_0)$ is gravity. In the numerical implementation, these dynamics are discretized using the multiple-shooting and RK4 procedure described in Section~\ref{ssec:discretization}.

The trajectory is required to satisfy the fixed boundary conditions
\begin{align*}
    x(0)=x_{\mathrm i},
    \qquad
    x(t_{\mathrm f})=x_{\mathrm f}.
\end{align*}

\paragraph*{Continuous-time path constraints}
The example includes three continuous-time path constraints, namely a tilt-angle bound, a thrust-magnitude bound, and a speed bound. These are encoded as STL predicates and imposed through an \always{} operator.

The tilt-angle constraint is represented by the cone condition
\vspace{-0.09cm}
\begin{align*}
    u_x^2 + u_y^2 \le (\cos \theta_{\max}\, u_z)^2,
\end{align*}
which ensures that the thrust vector remains within a cone of half-angle $\theta_{\max}$ about the vertical axis. The thrust and speed bounds are
\vspace{-0.09cm}
\begin{align*}
    \|u\|_2 \le T_{\max},
    \qquad
    \|v\|_2 \le v_{\max}.
\end{align*}
Accordingly, define the predicates
\vspace{-0.09cm}
\begin{align*}
    \mu_{\theta}
    &:=
    \big(
        (\cos \theta_{\max}\, u_z)^2 - u_x^2 - u_y^2 \ge 0
    \big), \\
    \mu_{T}
    &:=
    \big(
        T_{\max}^2 - u_x^2 - u_y^2 - u_z^2 \ge 0
    \big), \\
    \mu_{v}
    &:=
    \big(
        v_{\max}^2 - v_x^2 - v_y^2 - v_z^2 \ge 0
    \big).
\end{align*}
The continuous-time path-constraint formula is then
\begin{align*}
    \varphi_{\mathrm{path}}
    :=
    \bm G_{[0,t_{\mathrm f}]}
    \big(
        \mu_{\theta} \wedge \mu_{T} \wedge \mu_{v}
    \big).
\end{align*}

\paragraph*{Charging-station task}
Let $r_{\mathrm c}\in\mathbb{R}^3$ denote the center of the charging station and let $d_{\mathrm c}>0$ denote its radius. Define
\begin{align*}
    \mu_{\mathrm s}
    &:=
    \big(
        v_{\mathrm safe}^2 - v_x^2 - v_y^2 - v_z^2 \ge 0
    \big), \\
    \mu_{\mathrm c}
    &:=
    \big(
        d_{\mathrm c}^2 - \|r-r_{\mathrm c}\|_2^2 \ge 0
    \big).
\end{align*}
The task specification is
\vspace{-0.09cm}
\begin{align}
    \varphi_{\mathrm{task}}
    :=
    \mu_{\mathrm s}\,\bm U_{[0,t_{\mathrm f}]}\,\mu_{\mathrm c},
    \label{eq:num_until_spec}
\end{align}
which requires the vehicle speed to remain below the threshold until the charging station is reached.

For visualization, we also define the signed charging-station margin
\begin{align*}
    m_{\mathrm c}(t)
    :=
    d_{\mathrm c}-\|r(t)-r_{\mathrm c}\|_2,
\end{align*}
which is positive inside the charging station, zero on its boundary, and negative outside. Although the optimization uses the smooth predicate $\mu_{\mathrm c}$ above, the signed margin $m_{\mathrm c}$ provides a more interpretable geometric measure of station entry.

The overall CT-STL specification is
\begin{align*}
    \varphi := \varphi_{\mathrm{path}} \wedge \varphi_{\mathrm{task}}.
\end{align*}

Table~\ref{tab:num_params} summarizes the numerical parameters used in the simulation.


\begin{table}[t]
\centering
\renewcommand{\arraystretch}{1.08}
\setlength{\tabcolsep}{4.5pt}
\begin{tabular}{l|l||l|l}
\hline
Parameter & Value & Parameter & Value \\ \hline
$K$ & $5$ & $N$ & $10$ \\
$t_{\mathrm f}$ & $10.7649~\mathrm{s}$ & $g_0$ & $9.806~\mathrm{m/s^2}$ \\
$v_{\mathrm safe}$ & $5~\mathrm{m/s}$ & $v_{\max}$ & $10~\mathrm{m/s}$ \\
$\theta_{\max}$ & $45^\circ$ & $d_{\mathrm c}$ & $0.2~\mathrm{m}$ \\
$r_{\mathrm c}$ & $(2.5,\,-10,\,0)~\mathrm{m}$ & $c_i\in\bm c$ & $0.005$ \\
$x_{\mathrm i}$ & $(-10,\,-10,\,0)~\mathrm{m}$ & $x_{\mathrm f}$ & $(10,\,10,\,0)~\mathrm{m}$ \\
$T_{\max}$ & $1.75\,g_0~\mathrm{kg \, m/s^2}$ & $m$ & $1~\mathrm{kg}$  \\ \hline
\end{tabular}
\caption{Parameters for the quadrotor charging-station example.}
\label{tab:num_params}
\end{table}

\subsection{Penalized discretization and prox-convex for solution}\label{sec:num_method}

Let
\begin{align*}
    Z := (x_1,\dots,x_K,u_1,\dots,u_K)
\end{align*}
collect the nodal decision variables, and let $\mathcal Z$ denote the convex set induced by the fixed boundary conditions.

Since this example is posed as a feasibility problem, no performance cost is used. Instead, we solve an exact-penalty reformulation that penalizes the multiple-shooting defects and CT-STL violation:

\begin{align}
    \underset{Z\in\mathcal Z}{\mathrm{minimize}}
    \quad
    J_{\mathrm{nl}}(Z),
    \label{eq:num_penalized_problem}
\end{align}
where
\begin{align}
    J_{\mathrm{nl}}(Z)
    &:=
    w_{\mathrm{dyn}}
    \sum_{k=1}^{K-1}
    \left\|
        x_{k+1}-\Phi_{t_k}^{t_{k+1}}(x_k,u_k,u_{k+1})
    \right\|_1
    \nonumber\\
    &\quad
    +
    w_{\mathrm{stl}}
    \left|
        -\Gamma_{\bm c}^{\varphi}(X,U;\bar t_1)
    \right|_+ ,
    \label{eq:num_penalized}
\end{align}
$|a|_+ := \max(a,0)$, and $w_{\mathrm{dyn}}, w_{\mathrm{stl}} > 0$ are penalty weights.

We solve \eqref{eq:num_penalized_problem} using the prox-convex method~\cite{proxconvex}. At iteration $j+1$, given the current iterate $Z^j$, we solve the convex subproblem
\begin{align}
    \underset{Z\in\mathcal Z}{\mathrm{minimize}} \quad
    G(Z)
    &+
    H\!\left(
        C(Z^j)+\nabla C(Z^j)(Z-Z^j)
    \right) \nonumber \\
    &+
    \frac{w_{\rm ptr}^j}{2}\|Z-Z^j\|_2^2,
    \label{eq:num_prox_subproblem}
\end{align}
where $G$ encodes the convex terms, $H\circ C$ represents the penalized smooth nonconvex terms, and $w_{\rm ptr}^j>0$ is the adaptive proximal weight. 

In each prox-convex iteration, \texttt{JAX}~\cite{bradbury2018jax} is used to construct the nonlinear terms of the subproblem, including the discretized dynamics and the CT-STL specifications, as well as the corresponding first-order derivative information needed for linearization. Each convex subproblem is then modeled in \texttt{CVXPY}~\cite{diamond2016cvxpy} and solved by \texttt{QOCO}~\cite{chari2025qoco}, a conic convex solver. The initialization is chosen as a smooth trajectory interpolating between the initial state and the terminal state.

\subsection{Results}\label{sec:num_results}

\begin{figure}[t]
\centerline{\includegraphics[scale=0.50]{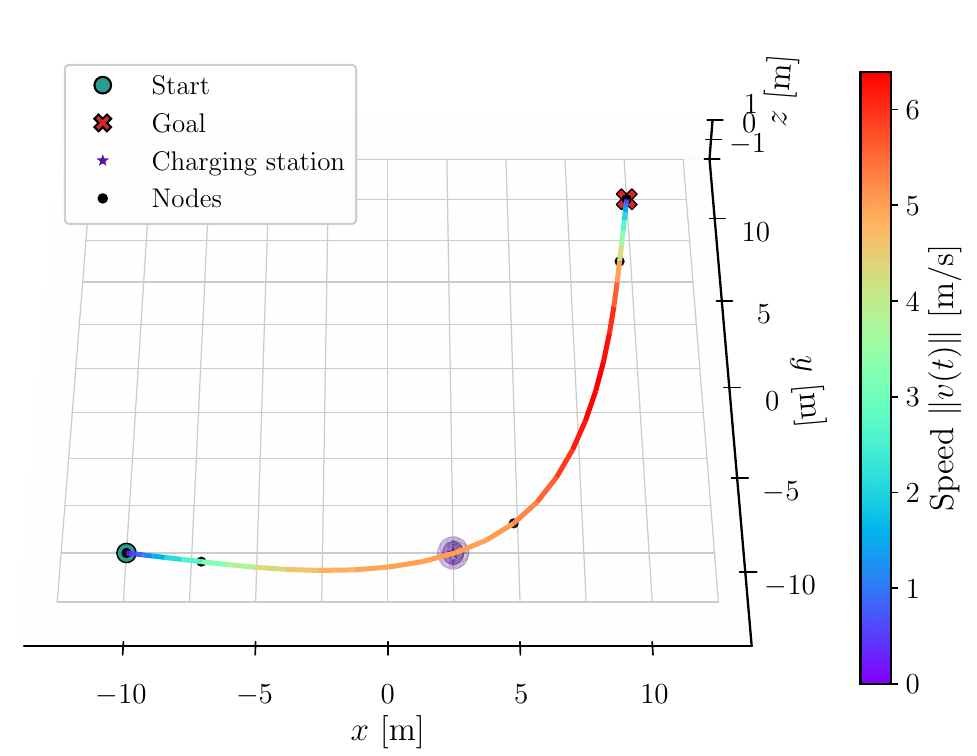}}
\caption{Optimized point-mass quadrotor trajectory together with the charging station. The vehicle reaches the charging region before proceeding to the terminal state, consistent with the continuous-time \until{} specification.}
\label{fig:num_traj}
\end{figure}

Figure~\ref{fig:num_traj} shows the optimized trajectory together with the charging station. The vehicle first approaches and enters the charging region, and then proceeds toward the terminal state. This behavior is consistent with the \until{} specification in \eqref{eq:num_until_spec}, since before the charging station is reached, the trajectory remains in the low-speed mode enforced by $\mu_{\mathrm s}$.

\begin{figure}[t]
\centerline{\includegraphics[scale=0.47]{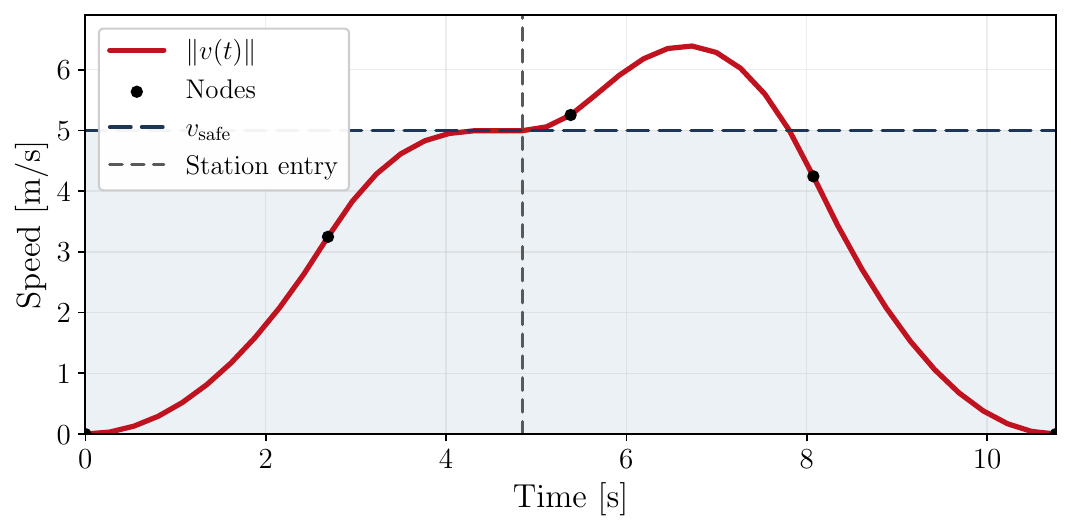}}
\caption{Speed profile along the optimized trajectory. The speed remains below the threshold $v_{\mathrm safe}$ until the charging station is reached, and always remains below the global bound $v_{\max}$.}
\label{fig:num_speed}
\end{figure}

Figure~\ref{fig:num_speed} shows the speed profile along the trajectory. The speed remains below the threshold $v_{\mathrm safe}$ prior to the charging event, as required by the \until{} semantics, while also satisfying the global speed bound encoded in $\varphi_{\mathrm{path}}$. After the charging station is reached, the trajectory is free to accelerate subject only to the path constraints.

\begin{figure}[t]
\centerline{\includegraphics[scale=0.47]{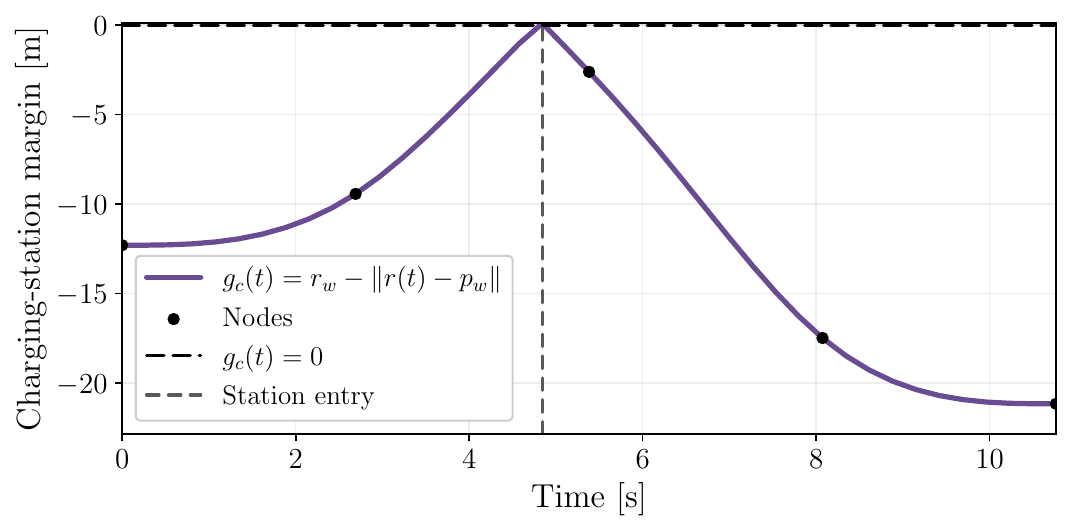}}
\caption{Signed charging-station margin $m_{\mathrm c}(t)=d_{\mathrm c}-\|r(t)-r_{\mathrm c}\|_2$. The margin is negative outside the station, zero on the boundary, and positive inside. The first zero crossing marks the witness time for the \until{} specification.}
\label{fig:num_margin}
\end{figure}

Figure~\ref{fig:num_margin} reports the signed charging-station margin $m_{\mathrm c}(t)=d_{\mathrm c}-\|r(t)-r_{\mathrm c}\|_2$. The margin is negative while the vehicle is outside the station, reaches zero at the station boundary, and becomes positive once the charging region is entered. The zero crossing therefore identifies the witness time at which the \until{} formula is satisfied.

The remaining \always{} constraints on thrust magnitude and tilt angle are also satisfied on the dense integration grid throughout the trajectory, although they are omitted from the plots for brevity. The final robustness value is positive, certifying satisfaction of the full CT-STL specification on the dense numerical trajectory representation.

The final iterate is dynamically consistent to numerical tolerance and satisfies the CT-STL specification on the dense integration grid. Using a simple, non-optimized Python implementation, the problem is solved in $10$ prox-convex iterations. The total time spent solving the convex subproblems is $15.35\,\mathrm{ms}$, and the total discretization time for the dynamics and CT-STL evaluations is $0.97\,\mathrm{ms}$, measured on a 2023 MacBook Pro with an Apple M2 processor. Because this implementation is intended only to demonstrate the method and has not been optimized for runtime, further gains from code optimization, compilation, and conversion to C/C++ are expected, with speedups on the order of $10\times$ to $40\times$.

\section{CONCLUSIONS AND FUTURE DIRECTIONS}

In this paper, we presented a smooth and exact parameterization of continuous-time Signal Temporal Logic (CT-STL) specifications for trajectory optimization. The proposed formulation evaluates temporal operators directly on the dense-time grid induced by the dynamics discretization, thereby enforcing CT-STL semantics on the same continuous-time trajectory representation used by the transcription. The resulting robustness measure is smooth, exact in sign, and well suited for gradient-based optimization, while alleviating the locality and masking issues of min/max-based smooth approximations.

As future work, we plan to develop a successive convexification framework for free-final-time optimal control problems with CT-STL specifications. In particular, we aim to incorporate a dynamic augmentation technique so that CT-STL satisfaction is tied directly to augmented dynamics and boundary conditions, avoiding the need to explicitly store the dense-time trajectory and allowing CT-STL constraints to be handled within a more standard direct optimal control framework. Dynamic augmentation has so far been available only for the \always{} operator in \cite{ctcs2024}; extending it to general CT-STL formulas is an important next step.


\bibliographystyle{ieeetr}      
\bibliography{ref}%

\end{document}

%% file: preamble.tex
\usepackage{xcolor}
\usepackage{amsmath} 
\usepackage{amssymb}
\usepackage{bm}
\usepackage{float}
\usepackage{booktabs}
\usepackage{mathtools}

\colorlet{fadedred}{red!30!gray!80!white}
\colorlet{fadedblue}{blue!70!gray!80!white}
\colorlet{fadedgreen}{green!60!blue!80!gray}
\usepackage{hyperref}
\hypersetup{
    colorlinks=true,
    linkcolor=fadedblue, 
    citecolor=fadedgreen,
    urlcolor=fadedgreen,
    pdftitle={ct-SCvx}
}

\newcommand{\always}{\texttt{always}}
\newcommand{\eventually}{\texttt{eventually}}
\newcommand{\implication}{\texttt{implication}}
\newcommand{\conjunction}{\texttt{conjunction}}
\newcommand{\disjunction}{\texttt{disjunction}}
\newcommand{\Negation}{\texttt{Negation}}

\newcommand{\until}{\texttt{until}}

\newcommand{\stlor}{\texttt{or}}
\newcommand{\stltrue}{\texttt{true}}
\newcommand{\stlfalse}{\texttt{false}}

\newcommand{\relu}[1]{|{#1}|_+}
\newcommand{\negrelu}[1]{|{#1}|_-}

\newcommand{\defeq}{\vcentcolon=}

\newtheorem{definition}{Definition}
\newtheorem{remark}{Remark}

\usepackage{tikz}
\usetikzlibrary{decorations.pathreplacing} 